\documentclass[reqno,11 pt]{amsart}
\usepackage{amsmath}
\usepackage{amssymb}
\input{diagrams}

\renewcommand {\SS}[1]{\mathfrak S_{#1}}

\newcommand {\ad}{\operatorname{ad}}
	
\newcommand {\ee}{\mathfrak e}
\newcommand {\so}{\mathfrak{so}}
\renewcommand {\sp}{\mathfrak{sp}}
\newcommand {\wt}{\widetilde}

\newcommand {\IC}{\mathbb{C}}
\newcommand {\IN}{\mathbb{N}}                          

\newcommand {\IY}{\mathbb Y}
\newcommand {\IZ}{\mathbb{Z}}
\newcommand {\g}{\mathfrak{g}}
\newcommand {\A}{\mathcal A}
\newcommand {\B}{\mathcal B}

\newcommand {\F}{\mathcal F}
\renewcommand {\H}{\mathcal H}
\newcommand {\M}{\mathcal M}
\newcommand {\V}{\mathcal V}

\newcommand {\half}[1]{\frac{#1}{2}}

\newcommand {\End}{\operatorname{End}}
\newcommand {\Ker}{\operatorname{Ker}}
\newcommand {\<}{\langle}
\renewcommand {\>}{\rangle}
\newtheorem {theorem}{Theorem}[section]
\newtheorem {proposition}[theorem]{Proposition}

\newtheorem {lemma}[theorem]{Lemma}
\renewcommand {\proof}{{\sc Proof.}\ }
\newcommand {\skproof}{{\sc Proof (sketch).}\ }

\newcommand {\remark}{\addtocounter{theorem}{1}{{\bf Remark \thetheorem.\ }}}
\newcommand {\definition}{\addtocounter{theorem}{1}{{\bf Definition \thetheorem.\ }}}
\newcommand {\example}{\addtocounter{theorem}{1}{{\bf Example \thetheorem.\ }}}

\newcommand {\KZ}{Knizhnik--Zamolodchikov  }
\newcommand {\KD}{Kohno--Drinfeld  }
\newcommand {\ie}{{\it i.e., }}
\newcommand {\halmos}{$\blacksquare$}

\newcommand {\Uh}{U_{\hbar}}
\newcommand {\Uhg}{\Uh\g}
\newcommand {\gl}[1]{\mathfrak{gl}_{#1}}
\renewcommand {\sl}[1]{\mathfrak{sl}_{#1}}

\newcommand {\Usl}[1]{U\sl{#1}}
\newcommand {\Uhgl}[1]{U_{\hbar}\gl{#1}}
\newcommand {\Uhsl}[1]{U_{\hbar}\sl{#1}}

\newcommand {\fml}{[\negthinspace[\hbar]\negthinspace]}

\newcommand {\ICh}{\IC\fml}
\newcommand {\Ch}{\IC_{\hbar}}
\renewcommand {\k}{^{(k)}}
\newcommand {\n}{^{(n)}}

\newcommand {\h}{\mathfrak h}
\newcommand {\reg}{_{\scriptscriptstyle{\operatorname{reg}}}}
\newcommand {\hreg}{\h\reg}
\newcommand {\Pg}{P_{\g}}
\newcommand {\Bg}{B_{\g}}

\newcommand {\Ug}{U\negthinspace\g}

\renewcommand {\dots}[1]{#1\cdots#1}

\newcommand {\tr}{\operatorname{tr}}

\newcommand {\eg}{{\it e.g., }}

\newcommand {\cvgt}{\{\negthinspace\negthinspace\{h\}\negthinspace\negthinspace\}}
\newcommand {\nfold}{V^{\otimes n}}

\newcommand{\cAl}{\circlearrowleft}
\newcommand{\car}{\circlearrowright}

\renewcommand {\a}{\mathfrak a}
\newcommand {\KKZ}{_{\scriptscriptstyle{\operatorname{KZ}}}}
\newcommand {\nablack}{\nabla_{\scriptscriptstyle{\operatorname{CKZ}}}}
\newcommand {\nablac}{\nabla_{\scriptscriptstyle{\operatorname{C}}}}
\newcommand {\nablak}{\nabla\KKZ}
\newcommand {\fmlo}{[\negthinspace[0]\negthinspace]}

\newcommand{\BMW}{\mathcal{BMW}}

\newcommand {\fd}{finite--dimensional }
\newcommand {\nth}{\negthickspace}

\newcommand {\Bn}{B_{n}}
\newcommand {\PW}{P_{W}}
\newcommand {\BW}{B_{W}}
\newcommand {\HW}{\H_{W}}
\newcommand {\wh}[1]{\widehat{#1}}
\newcommand {\Y}{\mathbb{Y}}
\newcommand {\nn}{\mathfrak n}
\newcommand {\Cg}{{\mathcal C}_{\g}}
\newcommand {\sign}{\operatorname{sign}}

\begin{document}

\title
{Flat Connections and Quantum Groups}
\author[V. Toledano Laredo]{Valerio Toledano Laredo}
\thanks{Work partially supported by an MSRI postdoctoral fellowship
for the academic year 2000--2001}
\begin{abstract}
We review the Kohno--Drinfeld theorem and a conjectural analogue
relating quantum Weyl groups to the monodromy of a flat connection
$\nablac$ on the Cartan subalgebra of a complex, semi-simple
Lie algebra $\g$ with poles on the root hyperplanes and values
in any $\g$-module $V$. We sketch our proof of this conjecture
when $\g=\sl{n}$ and when $\g$ is arbitrary and $V$ is a vector,
spin or adjoint representation. We also establish a precise link
between the connection $\nabla_{C}$ and Cherednik's generalisation
of the \KZ connection to finite reflection groups.
\end{abstract}
\address{
Institut de Mathematiques de Jussieu	\newline
UMR 7586, Case 191 			\newline
175 rue du Chevaleret, 75013 Paris	\newline
{\sf toledano@math.jussieu.fr}
}
\maketitle

\section{Introduction}

The aim of this paper is to discuss a principle which first
arose in the work of Kohno and Drinfeld and states, roughly
speaking, that quantum groups are natural receptacles for the
monodromy of certain integrable, first order PDE's. Quite how
general this principle is I do not know, but, as I will try
to show, it does extend beyond its original formulation.\\

The following diagram gives an overview of the paper
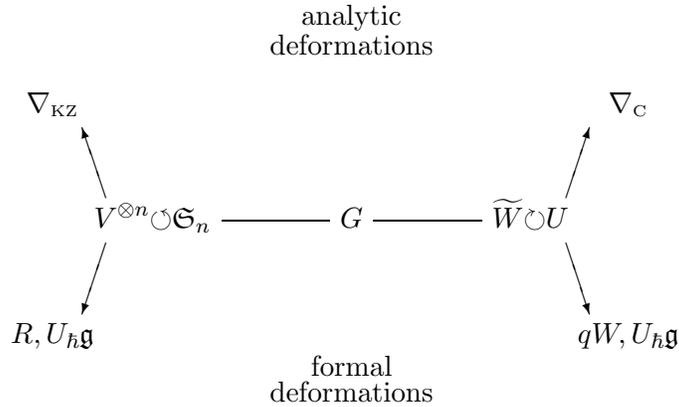
\newarrow{Ds}{-}{-}{-}{-}{-}
\newcommand {\anld}
{\stackrel{\textstyle{\text{analytic}}}{\text{deformations}}}
\newcommand {\fmld}
{\stackrel{\textstyle{\text{formal}}}{\text{deformations}}}
\begin{equation*}
\begin{diagram}[height=2em,width=.7em]
       &     &      &    &      &    &\anld&    &      &    & &     &       \\
\nablak&     &      &    &      &    &     &    &      &    & &     &\nablac\\
       &\luTo&      &    &      &    &     &    &      &    & &\ruTo&       \\
       &     &\nfold&\cAl&\SS{n}&\rDs&  G  &\rDs&\wt{W}&\car&U&     &       \\
       &\ldTo&      &    &      &    &     &    &      &    & &\rdTo&       \\
R,\Uhg &     &      &    &      &    &     &    &      &    & &     &qW,\Uhg\\
       &     &      &    &      &    &\fmld&    &      &    & &     &
\end{diagram}
\end{equation*}
\vskip .2cm

Here's how it should be read. To any complex, semi--simple
Lie group $G$ with Lie algebra $\g$, we may canonically
attach two finite groups. The first, the symmetric group
$\SS{n}$, is an external symmetry group and acts on the
$n$--fold tensor product of any finite--dimensional $G$--module
$V$. This action admits two remarkable deformations through
representations of Artin's braid group $\Bn$. The first is
the monodromy of the \KZ (KZ) equations, and is analytic in
the deformation parameter. The second is the $R$--matrix
representation of the Drinfeld--Jimbo quantum group $\Uhg$
associated to $\g$, and is formal. The remarkable theorem
of Kohno and Drinfeld alluded to above states that these two
seemingly very different deformations are in fact equivalent.\\

The second finite group attached to $G$ is its Weyl group $W$.
It is an internal symmetry group and it is tempting to think
that it acts on any finite--dimensional $G$--module $U$. This
isn't quite the case, but, as Tits showed \cite{Ti}, $W$
possesses a canonical abelian extension $\wt{W}$
\begin{equation*}
1\rightarrow\IZ_{2}^{r}\rightarrow\wt W\rightarrow W\rightarrow 1
\end{equation*}

by the sign group $\IZ_{2}^{r}$, with $r$ the rank of $G$,
which does act on $U$. This action is canonical only up to
conjugation by a fixed maximal torus $T$ of $G$, but since
this has little effect on the constructions I will discuss,
I will overlook this point and abusively speak of {\it the}
action of $\wt W$ on $U$.\\

Returning to the main story, this action possesses a formal
deformation through representations of the generalised braid
group $\Bg$ of type $\g$ known as the quantum Weyl group
action, which is constructed via the quantum group $\Uhg$.
It is natural to ask whether it also possesses an analytic
deformation obtained as the monodromy of a suitable flat
connection. As I will explain, the answer turns out to be
affirmative and is given by a new connection $\nablac$,
which I will call the Casimir connection. The latter was
discovered by De Concini around 1995 (unpublished), and
independently by J. Millson and myself \cite{MTL,TL2}, see
also \cite{FMTV}. The conjectural relation between these
two deformations, due to De Concini and myself will also be
discussed.\\

Here's a brief overview of the paper. In section \ref{se:hyperplane},
we describe a general method for constructing flat vector bundles
on hyperplane complements. This is applied in sections \ref{se:KZ}
and \ref{se:Casimir} to obtain the KZ and Casimir connections
respectively. Along the way, we describe in \S \ref{se:Coxeter}
Cherednik's generalisation of the KZ connection to other root
systems since it is closely related to the KZ and Casimir connections.
In section \ref{se:Uhg}, we briefly review the definition of the
quantum group $\Uhg$ and of the associated $R$--matrix and $q$Weyl
group representations. Sections \ref{se:KD} and \ref{se:DT} describe
the Kohno--Drinfeld theorem and its conjectural extension relating
the monodromy of the Casimir connection to $q$Weyl group actions.
We also sketch a proof of this conjecture for the case $\g=\sl{n}$,
referring to \cite{TL2} for more details. In section \ref{se:encore}
we study the relation of the Casimir and Cherednik connections.

\section{Flat connections on hyperplane complements}
\label{se:hyperplane}

Artin's braid groups and their generalised counterparts are, up to
the action of the corresponding finite Coxeter groups, fundamental
groups of {\it hyperplane complements}. This topological incarnation
leads to the analytic deformations mentioned in the Introduction
by taking the monodromy representations of suitable flat vector
bundles on these spaces. We begin by describing a general method
for constructing such bundles.\\

Recall that a hyperplane complement $X$ is defined by the following
data
\begin{itemize}
\item the base $\B$, a finite--dimensional complex vector space
\item the arrangement $\A=\{\H_{i}\}_{i\in I}$, a finite collection
of linear hyperplanes
\end{itemize}

and by setting $X=\B\setminus\A$. To describe flat vector bundles
over $X$, we need two additional pieces of data
\begin{itemize}
\item the fibre $\F$, a(nother) finite--dimensional complex vector
space
\item the residues $r_{i}\in\End(\F)$, labelled by the hyperplanes
in $\A$
\end{itemize}

With this at hand, we consider the following meromorphic connection
on the trivial vector bundle $\V=X\times\F$ over $X$
\begin{equation*}
\nabla=d-\sum_{i\in I}\frac{d\phi_{i}}{\phi_{i}}\cdot r_{i}
\end{equation*}

where $\phi_{i}\in\B^{*}$, $i\in I$, are linear equations for the
hyperplanes, so that $\H_{i}=\Ker(\phi_{i})$. The following useful
criterion of Kohno \cite{Ko1} tells us when such a connection is flat

\begin{lemma}\label{le:Kohno}
The above connection is flat iff, for any subcollection of linear
forms $\{\phi_{j}\}_{j\in J}$ which is maximal for the property
that their span in $\B^{*}$ is two--dimensional, one has
\begin{equation*}\label{eq:holo}
[r_{j},\sum_{j'\in J}r_{j'}]=0
\end{equation*}

for any $j\in J$.
\end{lemma}

The Lie theoretic nature of the above equations prompts
one to make, following Chen and Sullivan, the following\\

\definition The {\it holonomy Lie algebra} $\a(\A)$ of the
arrangement $\A$ is the quotient of the free Lie algebra
generated by symbols $r_{i}$, $i\in I$, by the relations
of lemma \ref{le:Kohno}.\\

Thus, if we decide to regard the $r_{i}$ as abstract
generators of $\a(\A)$ rather than endomorphisms of $\F$,
we may rephrase Kohno's lemma by saying that any linear
representation
\begin{equation*}
\pi:\a(\A)\longrightarrow\End(\F)
\end{equation*}

of $\a(\A)$ gives rise to a flat connection on $X\times\F$.
In fact, since the relations satisfied by the $r_{i}$ are
homogeneous, $\pi$ gives rise to a one--parameter family of
flat connections labelled by $h\in\IC$, namely
\begin{equation*}
\nabla=d-h\sum_{i\in I}\frac{d\phi_{i}}{\phi_{i}}\cdot r_{i}
\end{equation*}

and therefore to a one--parameter family of monodromy representations
of the fundamental group $\pi_{1}(X)$ of $X$. These analytically
deform the trivial representation of $\pi_{1}(X)$ on $\F$ which
is obtained by setting $h=0$. Thinking of this as a process of
exponentiation, we shall denote them by $e^{\pi}_{h}$ and think
of $\a(\A)$ as the Lie algebra of $\pi_{1}(X)$. Odd as this may
sound, this point of view is vindicated by the following

\begin{proposition}\label{pr:gen irred}
The monodromy representation $e^{\pi}_{h}:\pi_{1}(X)\longrightarrow
GL(\F)$ is generically irreducible iff the infinitesimal representation
$\pi:\a(\A)\longrightarrow\End(\F)$ is irreducible.
\end{proposition}

see, \eg \cite{MTL}. Here generically irreducible means irreducible
for all values of $h$ outside the zero set of some holomorphic
function $f\neq 0$.

\section{The Knizhnik--Zamolodchikov equations}\label{se:KZ}

Let $X_{n}$ be the configuration space of $n$ ordered points in $\IC$.
Thus
\begin{equation*}
X_{n}=\IC^{n}\setminus\bigcup_{1\leq i<j\leq n}\Delta_{ij}
\end{equation*}
where $\Delta_{ij}=\{(z_{1},\ldots,z_{n})\in\IC^{n}|z_{i}=z_{j}\}$
so that $X_{n}$ is a hyperplane complement. To construct the
\KZ (KZ) connection on $X_{n}$, we fix a complex, semi--simple
Lie algebra $\g$, one of its \fd representations $V$ and set
$\F=V^{\otimes n}$. The residue matrices $r_{ij}$ are usually
denoted by $\Omega_{ij}$ and are given by \cite{KZ}
\begin{equation*}
\Omega_{ij}=\sum_{a}\pi_{i}(X_{a})\pi_{j}(X^{a})
\end{equation*}

where $\{X_{a}\},\{X^{a}\}$ are dual basis of $\g$ with respect
to the {\it basic inner product} \ie the multiple $\<\cdot,\cdot
\>$ of the Killing form such that the highest root of $\g$ has
squared length $2$, and $\pi_{k}(X)$ denotes the action of $X
\in\g$ on the $k$th tensor factor in $V^{\otimes n}$. A simple
application of Kohno's lemma then shows that
\begin{equation*}
\nablak=d-
h\nth\sum_{\scriptscriptstyle{1\leq i<j\leq n}}\nth
\frac{d(z_{i}-z_{j})}{z_{i}-z_{j}}\cdot
\Omega_{ij}
\end{equation*}

is a flat connection on $X_{n}\times V^{\otimes n}$ for any $h\in
\IC$. Its monodromy yields a representation of Artin's pure braid
group on $n$ strands
\begin{equation*}
P_{n}=
\pi_{1}(\IC^{n}\setminus\{z_{i}=z_{j}\})
\longrightarrow
GL(V^{\otimes n})
\end{equation*}

which deforms the trivial representation of $P_{n}$ on $V^{\otimes
n}$. We can however do a little better by noticing that the symmetric
group $\SS{n}$ acts on $V^{\otimes n}$ and $X_{n}$. $\nablak$ is
readily seen to be equivariant for the combination of these two actions
and therefore descends to a flat connection on the quotient bundle
$\left(X_{n}\times V^{\otimes n}\right)/\SS{n}$ over $\wt X_{n}=
X_{n}/\SS{n}$ \ie the configuration space of $n$ unordered points
in $\IC$. Taking its monodromy, we obtain a one--parameter family
of representations of Artin's braid group on $n$ strands
\begin{equation*}
\rho_{h}:\Bn=\pi_{1}(\IC^{n}\setminus\{z_{i}=z_{j}\}/\SS{n})
\longrightarrow GL(V^{\otimes n})
\end{equation*}

$\rho_{h}$ depends analytically in $h$ and deforms the natural
action of $\SS{n}$ on $V^{\otimes n}$ since $\rho_{0}$ factors
through this action.\\

Recall that $\Bn$ is presented on elements $T_{i}$, $1\leq i\leq
n-1$ suject to Artin's braid relations \cite{Ar}
\begin{equation*}
\begin{array}{cl}
T_{i}T_{i+1}T_{i}=T_{i+1}T_{i}T_{i+1}& i=1\ldots n-1\\
T_{i}T_{j}=T_{j}T_{i}		     & |i-j|\geq 2
\end{array}
\end{equation*}

Each $T_{i}$ may be realised as a small loop in $\wt X_{n}$ around
the image of the hyperplane $\{z_{i}=z_{i+1}\}$. In particular, $
\rho_{h}(T_{i})$ is generically conjugate to $(i\medspace i+1)\cdot
\exp^{\pi\sqrt{-1}h\Omega_{i i+1}}$.\\

\example\label{ex:KZ gl} Take $\g=\gl{m}$ with vector representation
$V=\IC^{m}$ and basic inner product $\<X,Y\>=\tr_{V}(XY)$. If $e_{1},
\ldots,e_{n}$ is the standard basis of $V$ and $E_{ij}\medspace e_{k}
=\delta_{jk}e_{i}$ the corresponding elementary matrices then, on $V
^{\otimes 2}$,
\begin{equation*}
\Omega_{12}\medspace e_{k}\otimes e_{l}=
\sum_{1\leq i,j\leq m}
E_{ij}\otimes E_{ji}\medspace e_{k}\otimes e_{l}=
e_{l}\otimes e_{k}
\end{equation*}

so that $\Omega_{ij}$ acts on $V^{\otimes n}$ as the transposition
$(i\medspace j)$ and its eigenvalues are therefore $\pm 1$. The
corresponding monodromy representation
\begin{equation*}
\begin{diagram}[height=1.7em]
\Bn&     &\rTo^{\rho_{h}}&     &GL(V^{\otimes n}) \\
   &\rdTo&               &\ruTo&		  \\
   &     &\H_{\SS{n}}(q) &     & 
\end{diagram}
\end{equation*}

therefore factors through the Iwahori--Hecke algebra $\H_{\SS{n}}
(q)$, \ie the quotient of the group algebra $\IC[\Bn]$ by the
relations
\begin{equation*}
(T_{i}-q)(T_{i}+q^{-1})=0
\end{equation*}
where $q=e^{i\pi h}$.\\

\example Choose now an orthogonal vector space $V\cong\IC^{n}$, $\g=
\so(V)$ and let $e_{1},\ldots,e_{n}$ be an orthonormal basis of $V$.
Since the basic inner product on $\g$ is $\<X,Y\>=1/2\tr_{V}(XY)$,
we find
\begin{equation*}
\begin{split}
\Omega_{12}
&=
\sum_{1\leq i<j\leq n}
\left(E_{ij}-E_{ji}\right)\otimes\left(E_{ji}-E_{ij}\right)\\
&=
\sum_{1\leq i,j\leq n} E_{ij}\otimes E_{ji}-
\sum_{1\leq i,j\leq n} E_{ij}\otimes E_{ij}\\
&=
(1\medspace 2) - n\medspace p_{0}
\end{split}
\end{equation*}

where $p_{0}$ is the orthogonal projection onto the $\g$--fixed
line spanned by
\begin{equation*}
v_{0}=
\sum_{i=1}^{n}e_{i}\otimes e_{i}
\in S^{2}V
\end{equation*}

If, on the other hand, $V\cong\IC^{2n}$ is a sympletic vector space with
symplectic form $\omega$ and $\g=\sp(V)$, a similar computation in a
basis $e_{\pm 1},\ldots,e_{\pm n}$ of $V$ satisfying $\omega(e_{i},
e_{j})=\sign(i)\delta_{i+j,0}$ shows that 
\begin{equation*}
\Omega_{12}=(1\medspace 2) - 2n\medspace q_{0}
\end{equation*}
where $q_{0}$ is now the orthogonal projection onto the $\g$--fixed
line spanned by
\begin{equation*}
v_{0}=
\sum_{i=-n}^{n}\sign(i)\medspace e_{i}\otimes e_{-i}
\in\Lambda^{2}V
\end{equation*}

Thus, in either case, each generator of monodromy $T_{i}$ only has
the three eigenvalues $q,-q^{-1},r^{-1}$, where $q=e^{i\pi h}$ and
\begin{equation*}
r=\varepsilon e^{i\pi h(\dim(V)-\varepsilon)}
\qquad\text{with}\qquad
\varepsilon=
\left\{\begin{array}{rl}
+1&\text{if $V$ is orthogonal}\\
-1&\text{if $V$ is symplectic}
\end{array}\right.
\end{equation*}

With a little more work, one can show that the monodromy of $\nablak$
factors in this case through the Birman--Wenzl--Murakami algebra
$\BMW_{n}(q,r)$ \cite{BW,Mu} defined as the quotient of $\IC[\Bn]$
by the relations
\begin{gather*}
(T_{i}-q)(T_{i}+q^{-1})(T_{i}-r^{-1})=0\\
E_{i}T_{i-1}^{\pm 1}E_{i}=r^{\pm 1}E_{i}
\end{gather*}

where $E_{i}=1-(T_{i}-T_{i}^{-1})/(q-q^{-1})$ is a multiple of the
spectral projection of $T_{i}$ corresponding to the eigenvalue $r^{-1}$.


\section{The Coxeter--KZ connection}\label{se:Coxeter}

The connection described below was introduced by Cherednik \cite
{Ch}, to whom the results of this section are due, and is usually
referred to as the KZ connection. In order to distinguish it from
the one introduced in the previous section, we shall use the term
Coxeter--KZ (CKZ) connection instead. Let $W$ be a Weyl group, or
more generally a finite reflection group, with complexified reflection
representation $\h\cong\IC^{r}$ and root system $R=\{\alpha\}\subset
\h^{*}$. The base space and arrangement are now $\B=\h$ and
\begin{equation*}
\A=\bigcup_{\alpha\in R}\Ker(\alpha)
\end{equation*}
so that $X=\B\setminus\A$ is the space $\hreg$ of regular elements 
in $\h$. Set $\F=U$ where $U$ is a \fd $W$--module and let the
residue $r_{\alpha}$ be given by the reflection $s_{\alpha}\in W$.

\begin{theorem}
For any choice of weights $k_{\alpha}\in\IC$ satisfying
$k_{w\alpha}=k_{\alpha}$, $\forall w\in W$, the connection
\begin{equation*}
\nablack=
d-
\sum_{\alpha\succ 0}k_{\alpha}
\frac{d\alpha}{\alpha}\cdot s_{\alpha}
\end{equation*}
is a $W$--equivariant, flat connection on $\hreg\times U$.
\end{theorem}

The monodromy of $\nablack$ yields a family of representations
of the {\it generalised pure braid group $\PW$ of type $W$}
\begin{equation*}
\rho_{h}:\PW=\pi_{1}(\hreg)\longrightarrow GL(U)
\end{equation*}

deforming the trivial representation of $\PW$ on $U$. Each $W$--orbit
in $R$ carries a deformation parameter $k_{\alpha}$. As for the KZ
connection however, one can do a little better and use the action
of $W$ on $\hreg$ and $U$ to push $\nablack$ down to the quotient
$\hreg/W$. This yields a a representation of the {\it generalised
braid group of type $W$}
\begin{equation*}
\rho_{h}:\BW=\pi_{1}(\hreg/W)\longrightarrow GL(U)
\end{equation*}

which, for $k_{\alpha}=0$, factors through the action of $W$ on
$U$.\\

By Brieskorn's theorem \cite{Br}, $\BW$ is presented on generators
$S_{1},\ldots,S_{r}$ labelled by a choice of simple reflections
$s_{1},\ldots,s_{r}$ in $W$ with relations
\begin{equation*}
\underbrace{S_{i}S_{j}\cdots}_{m_{ij}}=
\underbrace{S_{j}S_{i}\cdots}_{m_{ij}}
\end{equation*}

for any $1\leq i<j\leq r$ where the number $m_{ij}$ of factors on
each side is equal to the order of $s_{i}s_{j}$ in $W$. Each $S_{i}$
may be obtained as a small loop in $\hreg/W$ around the reflecting
hyperplane $\Ker(\alpha_{i})$ of $s_{i}$ so that $\rho_{h}(S_{i})$
is generically conjugate to $s_{i}\exp^{\pi\sqrt{-1}k_{\alpha_{i}}
s_{i}}$. Since each simple reflection $s_{\alpha}$ has at most two
eigenvalues in $U$, the monodromy of $\nablack$
\begin{equation*}
\begin{diagram}[height=1.7em]
\BW&     &\rTo		&     &GL(U) \\
   &\rdTo&    		&\ruTo&      \\
   &     &\HW(q_{i})	&     &
\end{diagram}
\end{equation*}

therefore factors through the (unequal length) Hecke algebra $\HW
(q_{i})$ of $W$ \ie the quotient of $\IC[\BW]$ by the relations
\begin{equation*}\label{eq:unequal}
(S_{i}-q_{i})(S_{i}+q_{i}^{-1})=0
\end{equation*}

where $q_{i}=e^{\pi\sqrt{-1}\medspace k_{\alpha_{i}}}$. Choosing
$U$ to be the direct sum of the irreducible representations of $W$,
so that $\End(U)\cong\IC[W]$, and the weights $k_{\alpha}$ to be
generic, the monodromy does in fact yield an algebra isomorphism
of $\HW(q_{i})$ and $\IC[W]$.\\

\example When $W=\SS{n}$, the Coxeter--KZ connection is a particular
instance of the KZ connection. Indeed, we already noted in Example
\ref{se:KZ}.1. that, for $\g=\gl{m}$ acting on the $n$--fold tensor
product $V^{\otimes n}$ of its vector representation $V\cong\IC^{m}$,
the KZ operator $\Omega_{ij}$ is given by the transposition
$(i\medspace j)$. Thus,

\begin{proposition}\label{pr:KZ=CKZ}
The KZ connection for $\g=\gl{m}$ with values in $V^{\otimes n}$
coincides with the Coxeter--KZ connection for $W=\SS{n}$ with
values in $V^{\otimes n}$ and weights given by $k_{\alpha}=h$.
\end{proposition}

A finer version of this statement may of course be obtained using
Schur--Weyl duality. If $\lambda=(\lambda_{1},\ldots,\lambda_{m})
\in\IN^{m}$ is a Young diagram with at most $m$ rows and such that
$|\lambda|=\sum_{i}\lambda_{i}$ is equal to $n$, the irreducible
representation of $\gl{m}$ with highest weight $\lambda$ is a summand
in $V^{\otimes n}$. The corresponding multiplicity space $M_{\lambda}
^{n}$ is an irreducible representation of $\SS{n}$ and the KZ and
CKZ connections with values in $M_{\lambda}^{n}$ coincide.


\section{The Casimir connection}\label{se:Casimir}

We shall now use the Lie algebra $\g$ in a rather different way.
Fix a Cartan subalgebra $\h\subset\g$ and let $R=\{\alpha\}
\subset\h^{*}$ be the corresponding root system. The base
space and arrangement are the same as those of the Coxeter--KZ
connection for the Weyl group $W$ of $\g$, so that
\begin{equation*}
X=\h\setminus\bigcup_{\alpha\in R}\Ker(\alpha)=\hreg
\end{equation*}

The fibre $\F$ of the vector bundle is now a \fd $\g$--module
$U$. To describe the residue matrices $r_{\alpha}$, recall that
for any root $\alpha$, there is a corresponding subalgebra $\sl
{2}^{\alpha}\subseteq\g$ spanned by the triple $e_{\alpha},f_{
\alpha},h_{\alpha}$, where $h_{\alpha}=\alpha^{\vee}\in\h$ is the
corresponding coroot and $e_{\alpha},f_{\alpha}$ are a choice
of root vectors normalised by $[e_{\alpha},f_{\alpha}]=h_{\alpha}$.
The restriction of the basic inner product $\<\cdot,\cdot\>$ of
$\g$ to $\sl{2}^{\alpha}$ determines a canonical Casimir element
\begin{equation*}
C_{\alpha}=
\frac{\<\alpha,\alpha\>}{2}\left(
e_{\alpha}f_{\alpha}+f_{\alpha}e_{\alpha}+\half{1}h_{\alpha}^{2}
\right)
\in U\sl{2}^{\alpha}\subseteq \Ug
\end{equation*}

which we shall use as the residue on the hyperplane $\Ker(\alpha)$.
The following result was discovered by De Concini around 1995
(unpublished), and independently by J. Millson and myself \cite
{MTL,TL2}, see also \cite{FMTV}.

\begin{theorem} For any $h\in\IC$, the Casimir connection
\begin{equation*}
\nabla_{C}=
d-h\sum_{\alpha\succ 0}\frac{d\alpha}{\alpha}\cdot C_{\alpha}
\end{equation*}
is a flat connection on $\hreg\times U$ which is reducible
with respect to the weight space decomposition of $U$.
\end{theorem}
\proof Kohno's flatness criterion translates into the
statement that $\nabla_{C}$ is flat iff for any rank
2 root system $R_{2}\subseteq R$ determined by the
intersection of $R$ with a two--dimensional plane in
$\h^{*}$, the following holds for any positive root
$\alpha\in R_{2}$,
\begin{equation*}\label{eq:commutator}
[C_{\alpha},\sum_{\beta\in R_{2},\beta\succ 0}C_{\beta}]=0
\end{equation*}
Our original proof of this statement was a cumbersome
case--by--case check for the root systems $R_{2}=A_{1}
\times A_{1},A_{2},B_{2},G_{2}$. This was immediately
made obsolete by A. Knutson's elegant observation that
the second term in the commutator above is, modulo terms in
$\h$, the Casimir operator of the subalgebra $\g_{2}\subseteq
\g$ with root system $R_{2}$ and therefore commutes with $C_
{\alpha}$. The reducibility of $\nabla_{C}$ follows from the
fact that the $C_{\alpha}$ commute with $\h$ \halmos\\

We now wish to push the Casimir connection down to the quotient
$\hreg/W$ to get a monodromy representation of the generalised
braid group $\Bg=\BW$. This requires a little work because the
Weyl group $W$ does not act on $U$ and its Tits extension $\wt
{W}$, while acting on $U$, does not act freely on $\hreg$. To
circumvent this difficulty, we pull--back the Casimir connection
$\nablac$ to the universal cover $\wt{\hreg}\xrightarrow{p}
\hreg$. Since $\wt{W}$ is a quotient of $\Bg$, the latter acts
on $U$ and, freely, on $\wt{\hreg}$. The desired one--parameter
family $\rho_{h}$ of representations is obtained by taking the
monodromy of the flat vector bundle $(\wt{\hreg}\times U,p^{*}
\nablac)/\Bg$. It factors through the action of $\wt{W}$ on $U$
for $h=0$.\\

\example\label{ex:ad} Let $V=\g$ be the adjoint representation of
$\g$ so that the zero weight space $V[0]$ of $V$ is the Cartan
subalgebra $\h$ of $\g$. $V[0]$ is acted upon by the Casimirs
$C_{\alpha}$ as well as the Weyl group $W$ of $\g$ and, if
$t\in\h$
\begin{equation*}
\begin{split}
C_{\alpha}\medspace t
&= \<\alpha,\alpha\>\ad(e_{\alpha})\ad(f_{\alpha})\medspace t\\
&= \<\alpha,\alpha\>\<\alpha,t\> h_{\alpha}\\
&= \<\alpha,\alpha\>(1-s_{\alpha})t
\end{split}
\end{equation*}

From this we conclude that the Casimir connection with values
in $\h=V[0]$ coincides with the Coxeter--KZ connection with
values in the reflection representation of $W$, provided the
weights $k_{\alpha}$ are given by $-h\<\alpha,\alpha\>$ and
we tensor the CKZ connection with the character of $\pi_{1}
(\hreg/W)$ given by the multi--valuedness of the function
\begin{equation*}
f=\prod_{\alpha\succ 0}\alpha^{h\<\alpha,\alpha\>}
\end{equation*}

One cannot expect a similar coincidence to arise on the zero
weight space of any $\g$--module $V$ because the monodromy of
the CKZ connection with values in $V[0]$ always factors through
the Hecke algebra of $W$ while simple calculations show that
that of the Casimir connection hardly ever does. We shall
however return to this point in section \ref{se:encore}.\\

\remark Using the rigidity of the Hecke algebra of $W$ \ie
the fact that its representations are uniquely determined by
their specialisation at $q_{i}=1$ it is easy to see that the
monodromy representation of $\Bg$ on $\h=\g[0]$ is equivalent
to the reduced Burau representation of $\Bn=\Bg$ when $\g=\sl
{n}$ \cite{Bi} and to the Squier representation of $\Bg$ when
$\g$ is simply--laced \cite{Sq}.\\


%
%

\remark It is tempting to think that, since the Casimir operators
$C_{\alpha}$ are self--adjoint in any \fd $\g$--module, the connection
$\nablac$ is unitary whenever $h$ is purely imaginary. I am grateful
to P. Boalch for slapping my fingers on this point and pointing
out that this isn't (of course) so. Determining the values of $h$
for which the Casimir connection is unitary seems a very interesting
problem.

\section{Formal Deformations via Quantum Groups}\label{se:Uhg}

We turn now to formal deformations. These will be obtained via
the Drinfeld--Jimbo quantum group $\Uhg$. Recall that the latter
is a deformation of the enveloping algebra $\Ug$ of $\g$, \ie a
Hopf algebra over the ring $\ICh$ of formal power series in the
variable $\hbar$, which is topologically free as $\ICh$--module
and endowed with an isomorphism $\Uhg/\hbar\Uhg\cong\Ug$ of Hopf
algebras.\\

The simplest of these quantum groups corresponds to $\g=\sl{2}$,
where the standard generators $e,f,h$ of $\g$ given by 
\begin{equation*}
e=\begin{pmatrix}0&1\\0&0\end{pmatrix}\qquad
f=\begin{pmatrix}0&0\\1&0\end{pmatrix}\qquad
h=\begin{pmatrix}1&0\\0&-1\end{pmatrix}
\end{equation*}
together with the relations
\begin{gather*}
[h,e]=2e\qquad
[h,f]=-2f\\
[e,f]=h
\end{gather*}
which they satisfy are replaced by the generators $E,F,H$ of
$\Uhsl{2}$ subject to
\begin{gather*}
[H,E]=2E\qquad
[H,F]=-2F\\
\intertext{and}
[E,F]=
\frac{e^{\hbar H}-e^{-\hbar H}}{e^{\hbar}-e^{-\hbar}}
\end{gather*}

At first sight, the representation theory of $\Uhg$ offers few new
features. This is so because any \fd representation $\V$ of $\Uhg$,
\ie one which is finitely generated and topologically free as $\ICh
$--module, is uniquely determined by the $\g$--module $V=\V/\hbar\V$.
Indeed, since $H^{2}(\g,\Ug)=0$, the multiplication in $\Ug$ does
not possess non--trivial deformations and $\Uhg$ is isomorphic as
$\ICh$--algebra to
\begin{equation*}
\Ug\fml=\{\sum_{n\geq 0}x_{n}\hbar^{n}|\thinspace x_{n}\in\Ug\}
\end{equation*}
Using this to let $\Ug$ act on $\V$, we may regard the latter as a
deformation of $V$. Since $H^{1}(\g,End(V))=0$ however, $\V$ is
isomorphic, as $\Ug$ and therefore as $\Uhg$--module, to the trivial
deformation $V\fml$ of $V$.\\

The first novelty arises when one considers the action of the
symmetric group $\SS{n}$ on tensor products of $\g$--modules.
When implemented on the $n$--fold tensor product $\V^{\otimes n}$
of a \fd $\Uhg$--module $\V$, the latter does not commute with
the action of $\Uhg$. The following result shows however that
this problem may be corrected by deforming the action of $\SS{n}$.

\begin{theorem}[Faddeev--Reshetikhin--Takhtajan, Drinfeld, Jimbo]
There exists a universal $R$--matrix $R\in\Uhg\otimes\Uhg$ such
that the elements
$R_{i}^{\vee}\in GL(\V^{\otimes n})$, $i=1\ldots n-1$, given by
\begin{equation*}\label{eq:Ri}
R_{i}^{\vee}=
(i\medspace i+1)\cdot
\underbrace{1\otimes\cdots\otimes 1}_{i-1}
\otimes R\otimes
\underbrace{1\otimes\cdots\otimes 1}_{n-i-1}
\end{equation*}

commute with $\Uhg$ and satisfy
\begin{enumerate}
\item the braid relations :
\begin{equation*}
\begin{array}{cl}
R_{i}^{\vee}R_{i+1}^{\vee}R_{i}^{\vee}=
R_{i+1}^{\vee}R_{i}^{\vee}R_{i+1}^{\vee}
&i=1\ldots n-1\\[1.5 ex]
R_{i}^{\vee}R_{j}^{\vee}=R_{j}^{\vee}R_{i}^{\vee}
&|i-j|\geq 2
\end{array}
\end{equation*}
\item the deformation property : 
\begin{equation*}
R_{i}^{\vee}=(i\thickspace i+1)+o(\hbar)
\end{equation*}
\end{enumerate}
\end{theorem}

Thus, if one is prepared to replace $\SS{n}$ by the braid group
$\Bn$, there is an interesting, $\Uhg$--equivariant, 'permutation'
action of $\Bn$ on $\V^{\otimes n}$ which, when reduced mod $\hbar$,
factors through the natural action of $\SS{n}$. Moreover, this action
is {\it local} in the sense that $i$th generator of $\Bn$ acts on
the $i$ and $i+1$ tensor copies in $\V^{\otimes n}$ only, as does
the transposition $(i\medspace i+1)$. What is lost in this replacement
is the fact that the $R_{i}^{\vee}$ do not square to 1 and do not
therefore give an action of $\SS{n}$.\\


A similar phenomenon occurs for the action of the Tits extension
$\wt W$ on finite--dimensional $\g$--modules. There is no known,
canonical way to implement it on $\Uhg$--modules, but one may define
an action of the braid group $\Bg$ on these, known as the quantum
Weyl group action, which deforms that of $\wt W$. Before stating
the precise result, recall that the latter action arises by mapping
$\wt W$ to the completion $\wh{\Ug}$ of $\Ug$ with respect to its
finite--dimensional representations via
\begin{equation*}
s_{i}\longrightarrow\exp(e_{i})\exp(-f_{i})\exp(e_{i})
\end{equation*}

Let $q_{i}=e^{\hbar\<\alpha_{i},\alpha_{i}\>/2}$ and consider
the triple $q$--exponentials \cite{Ka,Sa}
\begin{gather*}\label{eq:triple q}
S_{i}=
\exp_{q_{i}^{-1}}(q_{i}^{-1}E_{i}q_{i}^{-H_{i}})
\exp_{q_{i}^{-1}}(-F_{i})
\exp_{q_{i}^{-1}}(q_{i}E_{i}q^{H_{i}})
\end{gather*}

where $E_{i},F_{i},H_{i}$ are the generators of the subalgebra
$\Uhsl{2}^{i}\subseteq\Uhg$ corresponding to the simple root
$\alpha_{i}$,
\begin{equation*}
\exp_{q}(x)=\sum_{n\geq 0}q^{n(n-1)/2}\frac{x^{n}}{[n]_{q}!}
\end{equation*}

and
\begin{equation*}
[n]_{q}=(q^{n}-q^{-n})/(q-q^{-1})
\qquad
[n]_{q}!=[n]_{q}[n-1]_{q}\cdots [1]_{q}
\end{equation*}

are the usual $q$--numbers and factorials. Viewing the $S_{i}$
as lying in the completion $\wh{\Uhg}$ of $\Uhg$ with respect
to its finite--dimensional representations, we have the following

\begin{theorem}[Lusztig, Kirillov--Reshetikhin, Soibelman]
The elements $S_{1},\ldots,S_{r}$ satisfy
\begin{enumerate}
\item the braid relations :
\begin{equation*}
\underbrace{S_{i}S_{j}S_{i}\cdots}_{m_{ij}}=
\underbrace{S_{j}S_{i}S_{j}\cdots}_{m_{ij}}
\end{equation*}
\item the deformation property :
\begin{equation*}
S_{i}=s_{i}+o(\hbar)
\end{equation*}
\end{enumerate}
\end{theorem}

The quantum Weyl group action is given by the $S_{i}$. Just as
the operators $R_{i}^{\vee}$, each $S_{i}$ is local in that it
lies in the completion $\wh{\Uhsl{2}^{i}}$ of $\Uhsl{2}^{i}$,
and does not square to 1.


\section{Monodromy theorems for Artin's braid groups}\label{se:KD}

Let us summarise what we have found so far for Artin's braid group
$\Bn$. We let as usual $\g$ be a complex, semi--simple Lie algebra, 
$\V$ a \fd representation of the quantum group $\Uhg$ and $V$ the
$\g$--module $\V/\hbar\V$.
\begin{equation*}
\begin{diagram}[height=1.6em]
   &               &&GL(V\dots{\otimes}V\fml)\\
   &\ruTo^{\nablak}&&                        \\
\Bn&               &&                        \\
   &\rdTo_{R}      &&                        \\
   &               &&GL(\V\dots{\otimes}\V)
\end{diagram}
\end{equation*}

On the one hand, $\Bn$ acts on $V^{\otimes n}$ via the monodromy
of the KZ equations. The latter depends analytically on the
deformation parameter $h\in\IC$ and can therefore be regarded as
an action of $\Bn$ on $V^{\otimes n}\cvgt$. Forgetting about
convergence, we regard $h$ as a formal variable, which we rename
$\hbar/2\pi i$, and consider the monodromy of $\nablak$ as an
action of $\Bn$ on $V^{\otimes n}\fml$. On the other hand, $\Bn$
acts on $\V^{\otimes n}$ via the $R$--matrix representation of
$\Uhg$. One now has the following beautiful

\begin{theorem}[Kohno, Drinfeld]
The monodromy representation of the KZ equations on $V^{\otimes
n}\fml$ is equivalent to the $R$--matrix action of $\Bn$ on $\V
^{\otimes n}$.
\end{theorem}

One may wonder whether the stated equivalence could be promoted
to an equality and proved by a direct calculation. There are
several reasons why this cannot be so.
\begin{itemize}
\item The monodromy representation depends upon a number of choices,
most notably that of a base point in the configuration space $\wt{X}
_{n}$. Thus, the upper row is an equivalence class of representations
rather than a single one.
\item This, in a sense, is also true of the $R$--matrix representation.
Indeed, to implement the latter on $V^{\otimes n}\fml$ rather than
on $\V^{\otimes n}$, one has to choose an algebra isomorphism $\phi:
\Uhg\rightarrow\Ug\fml$ to make $\Uhg$ act on $V\fml$. As mentioned,
such a $\phi$ exists, but is only unique up to  conjugation by an
element $a\in\Ug\fml$ of the form $1+o(\hbar)$.
\item This last objection partially disappears if one works modulo
$\hbar^{2}$ since in that case there is a preferred algebra isomorphism
\begin{equation*}
\Uhg/\hbar^{2}\Uhg\longrightarrow \Ug\otimes\ICh/(\hbar^{2})
\end{equation*}
obtained by lifting the given isomorphism $\Uhg/\hbar\Uhg\cong\Ug$.
Even then however one finds that the monodromy representation,
when computed mod $\hbar^{2}$ in a basis of horizontal sections
of $\nablak$ is not local, contrary to the $R$--matrix action.
\end{itemize}

The stated equivalence is in fact given by a rather explicit,
albeit cohomological expression (Drinfeld's twist) which is
not $\g$--equivariant \cite{Dr2,Dr3,Dr4}. Thus, the monodromy
and $R$--matrix pictures are complementary. The first is obtained
from the representation theory of $\g$, and is non--local, the
second is obtained from the representation theory of $\Uhg$
and is local.\\

Finally, we remark that when read from top to bottom, the
\KD theorem gives a concise description of the monodromy
of the KZ equations while, when read from bottom to top,
it is a sort of Riemann--Hilbert theorem since it asserts
that the $R$--matrix representation of $\Bn$ is the
monodromy of a flat connection on the trivial bundle
over the configuration space $X_{n}$.\\

Let us summarise the previous theorem as the following\\

{\bf Kohno--Drinfeld Principle.} If $\nabla$ is a flat
connection depending on a deformation parameter, there
exists a quantum group describing its monodromy.

\section{Monodromy theorems for generalised braid groups}\label{se:DT}

Turning now to the generalised braid group $\Bg$, we have a
similar diagram
\begin{equation*}
\begin{diagram}[height=1.4em]
   &               &&GL(V\fml)\\
   &\ruTo^{\nablac}&&         \\
\Bg&               &&         \\
   &\rdTo_{qW}     &&         \\
   &               &&GL(\V)
\end{diagram}
\end{equation*}

where the top row is the monodromy of the Casimir connection,
regarded as depending formally on the deformation parameter
$h$, here renamed $\hbar/2\pi i$, while the bottom one is the
quantum Weyl group action of $\Bg$ on the $\Uhg$--module $\V$.
In the light of the Kohno--Drinfeld principle, it seems natural
to make the following\\

{\bf Monodromy Conjecture.} The monodromy of the Casimir
connection with values in $V\fml$ is equivalent to the
quantum Weyl group action of $\Bg$ on $\V$.\\

This conjecture was formulated by De Concini in unpublished
work around 1995 and independently by myself in \cite{TL1,TL2}.
The difficulties in promoting its statement to a conjectural
equality are the same as for the Kohno--Drinfeld theorem. In
this case, the lack of locality of the monodromy representation
means that, even when computing mod $\hbar^{2}$, the image of
a small loop around the hyperplane $\Ker(\alpha_{i})$ does not
lie in the completion $\wh{\Usl{2}^{i}}$ of the $\sl{2}$--subalgebra
corresponding to the simple root $\alpha_{i}$.\\

A number of things can be proved in support of the above
conjecture, namely
\begin{itemize}
\item It is true for all representations of $\g=\sl{2}$
where $\Bg\cong\IZ$.
\item The spectra of the generators of $\Bg$ agree in both
representations.
\item It is true mod $\hbar^{2}$.
\end{itemize}

Moreover, one has the following
\begin{theorem}[\cite{TL1}] The monodromy conjecture holds
for the following pairs $(\g,V)$
\begin{itemize}
\item All fundamental representations of $\g=\sl{n}$.
\item Vector representation of $\g=\so_{n},\sp_{n}$.
\item Spin representation(s) of $\g=\so_{n}$.
\item Minuscule representations of $\g=\ee_{6},\ee_{7}$.
\item The 7--dimensional representation of $\g=\g_{2}$.
\item Adjoint representation of any $\g$.
\end{itemize}
\end{theorem}
\skproof All listed representations, except for the adjoint one,
have the property that their weight spaces are one--dimensional.
This makes it possible to compute the monodromy representation
explicitly, since, when restricted to the pure braid group $\Pg$,
it is just a sum of one--dimensional characters. One the other hand,
it is easy to deform these same $V$ to representations of $\Uhg$
explicitly \cite{Res}, and therefore to compute the corresponding
quantum Weyl group action using the triple $q$--exponentials that
define it. One finds in this case that the two representations
are conjugate by a diagonal matrix. The adjoint representation
of $\g$ requires a little more work. We first break $\g$ up as
$\nn\oplus\h$ where $\nn=\nn_{-}\oplus\nn_{+}$ is the direct sum of
the upper and lower nilpotent subalgebras, and $\h$ is the Cartan
algebra. Since both $\h$ and $\nn$ are preserved by the two actions,
it suffices to prove the monodromy conjecture for each piece. Since
the weight spaces of $\nn$ are one--dimensional, the corresponding
monodromy representation of $\Bg$ is readily computed. For the $q$Weyl
group action, one uses Lusztig's explicit deformation of the adjoint
representation \cite{Lu1}. The equivalence on $\nn$ is readily obtained
from this. For $\h$ we use the fact that both representations factor
through the Hecke algebra $\H_{W}(q_{i})$. This was shown in Example
\ref{se:Casimir}.2. for the monodromy representation and is a simple,
and old observation of Lusztig and Levendorskii--Soibelman for the
$q$Weyl group action. The equivalence is then obtained from the
rigidity of the Hecke algebra and the fact that both representations
are deformations of the reflection action of $W$ on $\h$ \halmos\\

\remark The above list of representations contains, for any simple
$\g$, at least one generator of the representation ring of $\g$, \ie
a $V$ such that any \fd irreducible $\g$--module is contained in a
tensor power $V^{\otimes n}$ of $V$. The monodromy conjecture
would therefore be proved if one could show that it holds for
$V_{1}\otimes V_{2}$ whenever it holds for each of the tensor
factors. This seems difficult.\\

For the case of $\g=\sl{n}$, we can say more.

\begin{theorem}[\cite{TL2}]\label{th:monodromy} The monodromy
conjecture holds for all representations of $\g=\sl{n}$.
\end{theorem}
\skproof The basic idea, summarised in the diagram below, is
to use the duality between $\gl{k}$ and $\gl{n}$ obtained
from their joint action on $k\times n$ matrices to reduce
the monodromy conjecture for $\sl{n}$ to the Kohno--Drinfeld
theorem for $\sl{k}$.
\newarrow{Corr}{<}{-}{-}{-}{>}
\newarrow{Dash}{<}{dash}{}{dash}{>}
\begin{equation*}
\begin{diagram}[height=2.5em,width=3em]
              &     &\IC[x_{11},\ldots,x_{kn}]&     &              \\
              &\ldTo&                         &\rdTo&              \\
\nablak,\sl{k}&     &                         &     &\nablac,\sl{n}\\
\dCorr^{KD}   &     &                         &     &\dDash        \\
R,\Uhsl{k}    &     &                         &     &qW,\Uhsl{n}   \\
              &\luTo&                         &\ruTo&              \\
              &     &\Ch[X_{11},\ldots,X_{kn}]&     &              
\end{diagram}
\end{equation*}

Let then $\A=\IC[x_{11},\ldots,x_{kn}]$ be the algebra of polynomial
functions on the space of $k\times n$ matrices. $\A$ is well--known
to be multiplicity--free, see \eg \cite[\S 132]{Zh}. Specifically,
if $\A^{d}\subset\A$ is the subspace of homogeneous polynomials of
degree $d\in\IN$, one has
\begin{equation*}
\A^{d}=
\bigoplus_{\substack{\lambda\in\IY_{\min(k,n)},\\ |\lambda|=d}}
V_{\lambda}\k\otimes V_{\lambda}\n
\end{equation*}

where $\Y_{p}$ is the set of Young diagrams $\lambda=(\lambda_{1}
,\ldots,\lambda_{p})\in\IN^{p}$ with at most $p$ rows, $|\lambda|
=\sum_{i}\lambda_{i}$ and $V_{\lambda}^{(p)}$ is the simple $\gl{p}
$--module with highest weight $\lambda$. If $k\geq n$, which we
henceforth assume, this allows one to identify the $\gl{n}$--weight
space $V_{\lambda}\n[\mu]$ corresponding to a weight $\mu=(\mu_{1},
\ldots,\mu_{n})\in\IN^{n}$ to the space $M_{\lambda}^{\mu}$ of highest
weight vectors of weight $\lambda$ for the diagonal $\gl{k}$--action
on
\begin{equation*}
\IC^{\mu_{1}}[x_{11},\ldots,x_{k1}]
\dots{\otimes}
\IC^{\mu_{n}}[x_{1n},\ldots,x_{kn}]
\end{equation*}

where $\IC^{\mu_{j}}[x_{1j},\ldots,x_{kj}]$ is the space of
polynomials in $x_{1j},\ldots,x_{kj}$ which are homogeneous
of degree $\mu_{j}$ . An explicit computation then proves
the following

\begin{proposition}
Under this identification, the Casimir connection $\nablac$ for
$\g=\sl{n}$ with values in $V_{\lambda}\n[\mu]$ coincides with
the KZ connection $\nablak$ for $\g'=\sl{k}$ with values in $M
_{\lambda}^{\mu}$.
\end{proposition}

Thus, the identification

\begin{equation*}
\bigoplus_{\nu\in\SS{n}\mu}V_{\lambda}\n[\nu]
\longrightarrow
\bigoplus_{\nu\in\SS{n}\mu}M_{\lambda}^{\nu}
\end{equation*}

is equivariant for the monodromy actions of $\Bn$ given by the
Casimir and KZ connections respectively.\\

Turning now to the $q$--setting, the algebra $\A$ possesses a
non--commutative, graded deformation $\A_{\hbar}$ over $\ICh$
on which both $\Uhgl{k}$ and $\Uhgl{n}$ act, and which is
multiplicity free. This allows as before to identify a weight
space $\V_{\lambda}\n[\mu]$ for $\Uhgl{n}$ with a corresponding
space $\M_{\lambda}^{\mu}$ of singular vectors for $\Uhgl{n}$.
An explicit, but a little more involved computation shows that

\begin{proposition}\label{pr:qW=R}
Under the identification
\begin{equation*}
\bigoplus_{\nu\in\SS{n}\mu}\V_{\lambda}\n[\nu]
\longrightarrow
\bigoplus_{\nu\in\SS{n}\mu}\M_{\lambda}^{\nu}
\end{equation*}
the ($\Uhgl{n}$--)$q$Weyl group action of $\Bn$ on the left--hand
side coincides with the $R$--matrix action for $\Uhgl{k}$ on the
right--hand side.
\end{proposition}

Proposition \ref{pr:qW=R} is the $q$--analogue of the simple fact
that the action of $\SS{n}$ on
\begin{equation*}
\IC[x_{11},\ldots,x_{kn}]
=
\IC[x_{11},\ldots,x_{k1}]
\dots{\otimes}
\IC[x_{1n},\ldots,x_{kn}]
\end{equation*}
obtained by permuting the columns of a $k\times n$ matrix is equal
to the one obtained by right multiplying the matrix by a permutation
matrix in $GL(n)$. The former action is the classical limit of the
$R$--matrix action of $\Uhgl{k}$, the latter that of the quantum
Weyl group action of $\Uhgl{n}$.\\

Putting together the above two propositions together with
the Kohno--Drinfeld theorem for $\sl{k}$, one obtains the
monodromy conjecture for $\sl{n}$ \halmos

\section{The Casimir and Coxeter--KZ connection (encore)}
\label{se:encore}

In this section, we pursue the study of the relations between
the Casimir connection for a Lie algebra $\g$ and the Coxeter--KZ
connection for its Weyl group $W$. The calculation of Example
\ref{se:Casimir}.2. for the adjoint representation of $\g$
generalises as follows.

\begin{proposition}\label{pr:V[[0]]}
For any \fd $\g$--module $V$, define
\begin{equation*}
V\fmlo=
\{v\in V[0]|\thinspace e_{\alpha}^{2}\medspace v=0,
\thickspace\forall\alpha\succ 0\}
\end{equation*}
Then,
\begin{enumerate}
\item $V\fmlo$ is invariant under $W$ and the $C_{\alpha}$.
\item On $V\fmlo$, one has
\begin{equation*}
C_{\alpha}=\<\alpha,\alpha\>(1-s_{\alpha})
\end{equation*}
so that the Casimir connection for $\g$ with values in $V\fmlo$
coincides with the Coxeter--KZ connection for $W$ with values
in $V\fmlo$ and weights given by $k_{\alpha}=-h\<\alpha,\alpha\>$.
\end{enumerate}
\end{proposition}
\proof The $W$--invariance of $V\fmlo$ is clear. Its $C_{\alpha}
$--invariance follows from (ii). Any $v\in V\fmlo$ may be written
as $v_{0}^{\alpha}+v_{2}^{\alpha}$ where $v^{\alpha}_{i}$ lies in
the zero weight space of the irreducible $\sl{2}^{\alpha}$--module
$V_{i}$ of dimension $i+1$. (ii) then follows from the fact that
$s_{\alpha}$ and $2/\<\alpha,\alpha\>\medspace C_{\alpha}$ act as
multiplication by $(-1)^{i/2}$ and $i(i+2)/2$ on $V_{i}[0]$ and
$V_{i}$ respectively \halmos\\

Note that, if $V$ is a {\it small} representation in the sense
of Broer and Reeder \cite{Bro,Re}, \ie is such that $2\alpha$
is not a weight of $V$ for any root $\alpha$, then $V\fmlo=V[0]$.
This is the case of the adjoint representation for example. In
general however, $V\fmlo$ can be a proper, non--zero subspace
of $V[0]$.\\

Proposition \ref{pr:V[[0]]} raises the question of whether
every irreducible representation $U$ of $W$ may be realised
inside some $V\fmlo$. On the positive side, we have the following.

\begin{proposition}\hfill
\begin{enumerate}
\item If $U_{\lambda}$ is the simple $\SS{n}$--module corresponding,
via the Schur--Weyl parametrisation, to the Young diagram $\lambda$,
then
\begin{equation*}
U_{\lambda}\cong
V_{\lambda^{t}}[0]=
V_{\lambda^{t}}\fmlo
\end{equation*}
where $V_{\lambda^{t}}$ is the irreducible representation of $\sl{n}
$ with highest weight given by the tranposed Young diagram $\lambda
^{t}$.
\item For any $\g$, the equality $\h=\g[0]$ induces an inclusion of
$W$--modules
\begin{equation*}
\Lambda^{i}\h\longrightarrow\Lambda^{i}\g\fmlo
\end{equation*}
\end{enumerate}
\end{proposition}
\proof (i) Let $V\cong\IC^{n}$ be the vector representation of $\sl
{n}$ so that the $V_{\lambda^{t}}$ span all irreducible summands of
$V^{\otimes n}$ as $\lambda$ varies over all partitions of $n$. A
simple inspection shows that $V^{\otimes n}$ is a small representation
so that $V^{\otimes n}\fmlo=V^{\otimes n}[0]$ \cite{Re}. The
isomorphism $U_{\lambda}\cong V_{\lambda^{t}}[0]$ is a simple
corollary of Schur--Weyl duality due to Kostant and Gutkin
\cite{Ks,Gu1,Gu2}. (ii) follows from an easy calculation \halmos\\

On the negative side however, one has the following

\begin{proposition}
There exist irreducible representations of the orthogonal Weyl
groups $B_{n}=W(\so_{2n+1})$, $n\geq 2$ and $D_{n}=W(\so_{2n})$,
$n\geq 4$ which are not contained in any $V\fmlo$.
\end{proposition}

It seems an interesting problem to determine, for any $\g$, the
Springer parameters of the irreducible representations of $W$
which arise inside some $V\fmlo$. A further motivation for this
question comes from the following simple corollary of proposition
\ref{pr:V[[0]]}. Let $\{U_{i}\}_{i\in I}$ be the isomorphism classes
of irreducible representations of $W$ which may be realised inside
some $V\fmlo$, and let $P_{0}\in\IC[W]$ be the corresponding central
projection onto $\bigoplus_{i\in I}\End(U_{i})$. Let $\Cg$ be the
{\it Casimir algebra} of $\g$, \ie the subalgebra
\begin{equation*}\label{eq:Cas}
\Cg=\<C_{\alpha}\>_{\alpha\succ 0}\subset\Ug
\end{equation*}

of the enveloping algebra of $\g$ generated by the $C_{\alpha}$.
Then,

\begin{proposition}
The assignement
\begin{equation*}
C_{\alpha}\rightarrow
\<\alpha,\alpha\>P_{0}(1-s_{\alpha})P_{0}
\end{equation*}
extends uniquely to a surjective, $W$--equivariant algebra
homomorphism of the Casimir algebra $\Cg$ of $\g$ onto the
subalgebra
\begin{equation*}
P_{0}\IC[W]P_{0}=\bigoplus_{i\in I}U_{i}\otimes U_{i}^{*}
\end{equation*}
of $\IC[W]$.
\end{proposition}


\begin{thebibliography}{ZZZ}


\bibitem{Ar} E. Artin, {\it Theory of Braids}, Ann. of
Math. {\bf 48} (1947), 101--26.

\bibitem{Bi} J. S. Birman, {\it Braids, Links, and Mapping
Class Groups}. Annals of Mathematics Studies, No. 82. Princeton
University Press, Princeton, 1974.

\bibitem{BW} J. S. Birman, H. Wenzl, {\it Braids, Link
Polynomials and a New Algebra}, Trans. Amer. Math. Soc.
{\bf 313} (1989), 249--273. 


\bibitem{Br} E. Brieskorn, {\it Die Fundamentalgruppe
des Raumes der Regul\"aren Orbits einer Endlichen Komplexen
Spiegelungsgruppe}, Invent. Math. {\bf 12} (1971), 57--61.

\bibitem{Bro} A. Broer, {\it The Sum of Generalized
Exponents and Chevalley's Restriction Theorem for Modules
of Covariants}, Indag. Math. (N.S.) {\bf 6} (1995), 385--396.

\bibitem{Ch} I. V. Cherednik, {\it Generalized Braid
Groups and Local $r$--Matrix Systems}, Soviet Math. Dokl.
{\bf 40} (1990), 43--48.


\bibitem{Dr1} V. G. Drinfeld, {\it Quantum Groups},
Proceedings of the International Congress of Mathematicians,
Berkeley 1986, 798--820.

\bibitem{Dr2} V. G. Drinfeld, {\it On Almost Cocommutative
Hopf Algebras}, Leningrad Math. J. {\bf 1} (1990), 321--342. 

\bibitem{Dr3} V. G. Drinfeld, {\it Quasi--Hopf Algebras},
Leningrad Math. J. {\bf 1} (1990), 1419--57.

\bibitem{Dr4} V. G. Drinfeld, {\it On Quasitriangular
Quasi--Hopf Algebras and on a Group that is closely connected
with ${\rm Gal}(\overline{Q}/{Q})$}, Leningrad Math. J. {\bf 2}
(1991), 829--860.


\bibitem{FMTV} G. Felder, Y. Markov, V. Tarasov, A.
Varchenko, {\it Differential Equations Compatible with KZ
Equations}, Math. Phys. Anal. Geom. {\bf 3} (2000), 139--177.

\bibitem{FRT} L. D. Faddeev, N. Y. Reshetikhin, L. A. Takhtajan,
{\it Quantization of Lie groups and Lie algebras}, Algebraic analysis,
Vol. I, 129--139, Academic Press, Boston, MA, 1988.


\bibitem{Gu1} E. A. Gutkin, {\it Representations of the
Weyl Group in the Space of Vectors of Zero Weight}, Uspehi
Mat. Nauk {\bf 28} (1973), 237--238.

\bibitem{Gu2} E. A. Gutkin, {\it Schur--Weyl duality and
Representations of Permutation Groups}, preprint, October 2001.

\bibitem{Ho} R. Howe, {\it Remarks on Classical Invariant
Theory}. Trans. Amer. Math. Soc. {\bf 313} (1989), 539--570.

\bibitem{Ka} M. Kashiwara, {\it Private Note on
Finite--Dimensional Representations of Quantized
Affine Algebras}, unpublished notes.

\bibitem{KR} A. N. Kirillov, N. Reshetikhin,
{\it $q$-Weyl Group and a Multiplicative Formula
for Universal $R$-Matrices}, Comm. Math. Phys.
{\bf 134} (1990), 421--431.

\bibitem{KZ} V. G. Knizhnik, A. B. Zamolodchikov, {\it
Current Algebra and Wess--Zumino Model in Two Dimensions}.
Nuclear Phys. B {\bf 247} (1984), 83--103.

\bibitem{Ko1} T. Kohno, {\it Quantized Enveloping
Algebras and Monodromy of Braid Groups}, preprint, 1988.

\bibitem{Ko2} T. Kohno, {\it Integrable Connections
Related to Manin and Schechtman's Higher Braid Groups},
Illinois J. Math. {\bf 34} (1990), 476--484.

\bibitem{Ks} B. Kostant, {\it On Macdonald's $\eta
$--Function Formula, the Laplacian and Generalized Exponents},
Advances in Math. {\bf 20} (1976), 179--212. 


\bibitem{Lu1} G. Lusztig, {\it Finite--Dimensional
Hopf Algebras Arising from Quantized Universal Enveloping
Algebra}. J. Amer. Math. Soc. {\bf 3} (1990), 257--296.

\bibitem{Lu2} G. Lusztig, {\it Introduction to Quantum
Groups}. Progress in Mathematics, 110. Birkh\"auser Boston,
1993.


\bibitem{MTL} J. Millson, V. Toledano Laredo,
{\it Casimir Operators and Monodromy Representations
of Generalised Braid Groups}, in preparation.

\bibitem{Mu} J. Murakami, {\it The Kauffman Polynomial
of Links and Representation Theory}, Osaka J. Math. {\bf 24}
(1987), 745--758. 

\bibitem{Re} M. Reeder, {\it Zero Weight Spaces and
the Springer Correspondence}, Indag. Math. (N.S.) {\bf
9} (1998), 431--441.

\bibitem{Res} N. Reshetikhin, {\it Quantized Universal
Enveloping Algebras, the Yang-Baxter Equation and Invariants
of Links I}, LOMI preprint E--4--87, 1987.

\bibitem{Sa} Y. Saito, {\it PBW Basis of Quantized
Universal Enveloping Algebras}, Publ. Res. Inst. Math.
Sci. {\bf 30} (1994), 209--232.

\bibitem{So} Y. S. Soibelman, {\it Algebra
of Functions on a Compact Quantum Group and its
Representations}, Leningrad Math. J. {\bf 2}
(1991), 161--178.

\bibitem{Sq} C. C. Squier, {\it Matrix Representations
of Artin Groups}, Proc. Amer. Math. Soc. {\bf 103} (1988),
49--53.

\bibitem{TL1} V. Toledano Laredo, {\it Monodromy
Representations of Generalised Braid Groups and Quantum
Weyl Groups}, in preparation.

\bibitem{TL2} V. Toledano Laredo, {\it A Kohno--Drinfeld
Theorem for Quantum Weyl Groups}, Duke Math. J. {\bf 112}
(2002), 421--51.

\bibitem{Ti} J. Tits, {\it Normalisateurs de Tores.
I. Groupes de Coxeter Etendus}, J. Algebra {\bf 4}
(1966), 96--116.


\bibitem{Zh} D. P. Zhelobenko, {\it Compact Lie Groups and
their Representations.}, Translations of Mathematical
Monographs, Vol. 40. American Mathematical Society, 1973.

\end{thebibliography}
\end{document}